\documentclass{amsart}

\usepackage{amssymb}

\newtheorem{theorem}{Theorem}[section]
\newtheorem{lemma}[theorem]{Lemma}
\newtheorem{proposition}[theorem]{Proposition}
\newtheorem{corollary}[theorem]{Corollary}

\theoremstyle{definition}
\newtheorem{definition}[theorem]{Definition}
\newtheorem{example}[theorem]{Example}

\theoremstyle{remark}
\newtheorem{remark}[theorem]{Remark}

\numberwithin{equation}{section}

\newcommand{\bdf}{\begin{definition}}
\newcommand{\edf}{\end{definition}}

\newcommand{\blem}{\begin{lemma}}
\newcommand{\elem}{\end{lemma}}

\newcommand{\bthm}{\begin{theorem}}
\newcommand{\ethm}{\end{theorem}}

\newcommand{\bpf}{\begin{proof}}
\newcommand{\epf}{\end{proof}}

\newcommand{\bprop}{\begin{proposition}}
\newcommand{\eprop}{\end{proposition}}

\newcommand{\bcor}{\begin{corollary}}
\newcommand{\ecor}{\end{corollary}}

\newcommand{\brem}{\begin{remark}}
\newcommand{\erem}{\end{remark}}

\newcommand{\bquest}{\begin{question}}
\newcommand{\equest}{\end{question}}

\newcommand{\bex}{\begin{example}}
\newcommand{\eex}{\end{example}}

\newcommand{\benu}{\begin{enumerate}\renewcommand{\labelenumi}{{\rm (\arabic{enumi})}}\renewcommand{\itemsep}{0pt}}
\newcommand{\eenu}{\end{enumerate}}

\newcommand{\hil}{\mathcal{H}}

\newcommand{\C}{\mathbb{C}}
\newcommand{\free}{{\mathbb{F}_d}}
\newcommand{\lp}{{\ell_p}}

\begin{document}

\title{Free group $C^*$-algebras associated with $\ell_p$}

\author{Rui OKAYASU}
\address{Department of Mathematics Education, Osaka Kyoiku University,
Kashiwara, Osaka 582-8582, JAPAN}
\email{rui@cc.osaka-kyoiku.ac.jp}
\thanks{The author was supported in part by JSPS}


\subjclass[2000]{Primary 46L05; Secondary 22D25}

\date{}


\keywords{$C^*$-algebras, free groups, positive definite functions}

\begin{abstract}
For every $p\geq 2$, we give a characterization of positive definite functions on a free group with finitely many generators, which can be extended to the positive linear functionals on the free group $C^*$-algebra associated with the ideal $\ell_p$. This is a generalization of Haagerup's characterization for the case of the reduced free group $C^*$-algebra. As a consequence, the associated $C^*$-algebras are mutually non-isomorphic, and they have a unique tracial state.
\end{abstract}

\maketitle

\section{Introduction}
N. P. Brown and E. Guentner introduce new $C^*$-completion of the group ring of a countable discrete group $\Gamma$ in \cite{bg}. More precisely, for a given algebraic two-sided ideal in $\ell_\infty(\Gamma)$, they define the associated group $C^*$-algebra. These recover the full group $C^*$-algebra for $\ell_\infty(\Gamma)$ itself, and the reduced group $C^*$-algebra for $c_c(\Gamma)$, respectively. Hence if we take $c_0(\Gamma)$ or $\ell_p(\Gamma)$ with $p\in[1,\infty)$ for example, we may obtain new group $C^*$-algebra. We remark that a standard characterization of amenability implies that the associated $C^*$-algebras of an amenable group are all isomorphic for any ideals. In \cite{bg}, they also give a characterization of the Haagerup property and Property (T) in terms of ideal completions. 

In this paper, we study their $C^*$-algebra of a free group associated with $\ell_p$. By \cite{bg}, for any $p\in[1,2]$, the group $C^*$-algebra of $\Gamma$ associated with $\ell_p$ is isomorphic to the reduced group $C^*$-algebra. Therefore the case where $p\in(2,\infty)$ is essential. Our main result is a characterization of positive definite functions on a free group, which can be extended to the positive linear functionals on the free group $C^*$-algebra associated with $\ell_p$. In \cite{ha}, U. Haagerup gives its characterization for the case of the reduced free group $C^*$-algebra. Thus our theorem is a generalization of his result. As a consequence, the free group $C^*$-algebras associated with $\ell_p$ are mutually non-isomorphic for every $p\in[2,\infty]$. We also obtain that for each $p\in[2,\infty)$, the free group $C^*$-algebra associated with $\ell_p$ has a unique tracial state. Moreover, we consider algebraic ideals
$$D_p^-(\Gamma)=\bigcup_{\varepsilon>0}\ell_{p-\varepsilon}(\Gamma)$$
for $1<p\leq\infty$, and
$$D_p^+(\Gamma)=\bigcap_{\varepsilon>0}\ell_{p+\varepsilon}(\Gamma)$$
for $1\leq p<\infty$. Then, by using our characterization, it is also shown that the free group $C^*$-algebra associated with $D_p^\pm$ coincides with the free group $C^*$-algebra associated with $\ell_p$. 

\section{Preliminaries}\label{pre}
In this section, we fix the notations for the convenience of the reader and recall some results in \cite{bg}. 

Let $\Gamma$ be a countable discrete group and $\pi$ be a unitary representation of $\Gamma$ on a Hilbert space $\hil$. For $\xi,\eta\in\hil$, we denote the matrix coefficient of $\pi$ by
$$\pi_{\xi,\eta}(s)=\langle \pi(s)\xi,\eta\rangle.$$
Note that $\pi_{\xi,\eta}\in\ell_\infty(\Gamma)$, where $\ell_\infty(\Gamma)$ is the abelian $C^*$-algebra of all bounded functions on $\Gamma$. 

Let $D$ be an algebraic two-sided ideal of $\ell_\infty(\Gamma)$. If there exists a dense subspace $\hil_0$ of $\hil$ such that $\pi_{\xi,\eta}\in D$ for all $\xi, \eta\in\hil_0$, then $\pi$ is called $D$-representation. If $D$ is invariant under the left and right translation of $\Gamma$ on $\ell_\infty(\Gamma)$, then it is said to be translation invariant. 

Throughout this paper, we assume that $D$ is a non-zero translation invariant ideal of $\ell_\infty(\Gamma)$. For each $p\in[1, \infty)$, we denote  the norm on $\ell_p(\Gamma)$ by 
$$|f|_p=\left(\sum_{s\in\Gamma}|f(s)|^p\right)^{\frac{1}{p}}\ \ \ \mbox{for}\ f\in\ell_p(\Gamma).$$
Note that $\ell_p(\Gamma)$ is a translation invariant ideal of $\ell_\infty(\Gamma)$. We denote by $c_0(\Gamma)$ the of functions on $\Gamma$, vanishing at infinity. It is the non-trivial closed translation invariant ideal of $\ell_\infty(\Gamma)$. 

Under our assumption, $D$ contains $c_c(\Gamma)$, which is the ideal of all finitely supported functions on $\Gamma$. Moreover, if $\pi$ has a cyclic vector $\xi$ such that $\pi_{\xi,\xi}\in D$, then $\pi$ is a $D$-representation with respect to a dense subspace 
$$\hil_0=\mathrm{span}\{\pi(s)\xi : s\in\Gamma\}.$$ 
We denote by $\lambda$ the left regular representation of $\Gamma$. It is easy to see that $\lambda$ is a $c_c$-representation, or a $D$-representation for any $D$. 

The $C^*$-algebra $C_D^*(\Gamma)$ is the $C^*$-completion of the group ring $\mathbb{C}\Gamma$ by $\|\cdot\|_D$, where
$$\|f\|_D=\sup\{\|\pi(f)\| : \pi\ \mbox{is a $D$-representation}\}\ \ \ \mbox{for}\ f\in\C\Gamma.$$
Note that if $D_1$ and $D_2$ are ideals of $\ell_\infty(\Gamma)$ with $D_1\supset D_2$, then there exists the canonical quotient map from $C_{D_1}^*(\Gamma)$ onto $C_{D_2}^*(\Gamma)$. We denote by $C^*(\Gamma)$ the full group $C^*$-algebra, and by $C_\lambda^*(\Gamma)$ the reduced group $C^*$-algebra, respectively. In \cite{bg}, the following results are obtained:

\begin{itemize}
\item $C^*(\Gamma)= C_{\ell_\infty}^*(\Gamma)$ and $C_\lambda^*(\Gamma)= C_{c_c}^*(\Gamma)$. \\
\item $C_{\ell_p}^*(\Gamma)= C_\lambda^*(\Gamma)$ for every $p\in[1,2]$. \\
\item $C^*(\Gamma)= C_D^*(\Gamma)$ if and only if there exists a sequence $(h_n)$ of positive definite functions in $D$ such that $h_n\to 1$. \\
\item If $C^*(\Gamma)= C_{\ell_p}^*(\Gamma)$ for some $p\in[1, \infty)$, then $\Gamma$ is amenable. \\
\item $\Gamma$ has the Haagerup property if and only if $C^*(\Gamma)= C_{c_0}^*(\Gamma)$.
\end{itemize}

\section{Positive definite functions on a free group}

Let $\free$ be the free group on finitely many generators $a_1,\dots,a_d$ with $d\geq 2$. We denote by $|s|$ the word length of $s\in\free$ with respect to the canonical generating set $\{a_1, a_1^{-1},\dots, a_d, a_d^{-1}\}$. For $k\geq 0$, we put 
$$W_k=\{s\in\free : |s|=k\}.$$ 
We denote by $\chi_k$ the characteristic function for $W_k$. 

In the following lemma, the case where $q=2$ is given by Haagerup in {\cite[Lemma 1.3]{ha}}. His proof also works for $q\in[1,2]$ by using H\"older's inequality, instead of Cauchy-Schwarz inequality. We remark that this is also appeared in \cite{bo}. 

\blem\label{conv}
Let $q\in[1,2]$. Let $k,\ell$ and $m$ be non negative integers. Let $f$ and $g$ be functions on $\free$ such that $\mathrm{supp}(f)\subset W_k$ and $\mathrm{supp}(g)\subset W_\ell$, respectively. If $|k-\ell|\leq m\leq k+\ell$ and $k+\ell-m$ is even, then 
$$|(f*g)\chi_m|_q\leq|f|_q|g|_q,$$
and if $m$ is any other value, then 
$$|(f*g)\chi_m|_q=0.$$
\elem

\bpf
It is shown by an argument similar as in {\cite[Lemma 1.3]{ha}}. However for convenience, we give the complete proof.

Note that
$$(f*g)(s)=\sum_{\substack{t,u\in\free \\ s=tu}}f(t)g(u)=\sum_{\substack{|t|=k \\ |u|=\ell \\ s=tu }}f(t)g(u).$$
Since the possible values of $|tu|$ are $|k-\ell|, |k-\ell|+2,\dots,k+\ell$, we have 
$$|(f*g)\chi_m|_q=0$$
for any other values of $m$.

The case where $q=1$ is trivial. So we consider the case where $q\ne 1$.

First we assume that $m=k+\ell$. In this case, if $|s|=m$, then $s$ can be uniquely written as a product $tu$ with $|t|=k$ and $|u|=\ell$. Hence
$$(f*g)(s)=f(t)g(u).$$
Therefore
$$|(f*g)\chi_m|_q^q=\sum_{\substack{|tu|=k+\ell \\ |t|=k \\ |u|=\ell}}|f(t)|^q|g(u)|^q\leq\sum_{\substack{|t|=k \\ |u|=\ell}}|f(t)|^q|g(u)|^q=|f|_q^q|g|_q^q.$$

Next we assume that $m=|k-\ell|, |k-\ell|+2,\dots,k+\ell-2$. In these cases, we have $m=k+\ell-2j$ for $1\leq j\leq\min\{k, \ell\}$. Let $s=tu$ with $|s|=m$, $|t|=k$ and $|u|=\ell$. Then $s$ can be uniquely written as a product $t'u'$ such that $t=t'v$, $u=v^{-1}u'$ with $|t'|=k-j$, $|u'|=\ell-j$ and $|v|=|v^{-1}|=j$. We define 
$$f'(t)=\left(\sum_{|v|=j}|f(tv)|^q\right)^{\frac{1}{q}} \ \mbox{if}\ |t|=k-j,\ \mbox{and}\ f'(t)=0\ \mbox{otherwise.}$$
We also define 
$$g'(u)=\left(\sum_{|v|=j}|g(v^{-1}u)|^q\right)^{\frac{1}{q}} \ \mbox{if}\ |u|=\ell-j,\ \mbox{and}\ g'(u)=0\ \mbox{otherwise.}$$
Note that $\mathrm{supp}(f')\subset W_{k-j}$ and $\mathrm{supp}(g')\subset W_{\ell-j}$. Moreover 
$$|f'|_q^q=\sum_{|t|=k-j}\left(\sum_{|v|=j}|f(tv)|^q\right)=|f|_q^q,$$
and similarly $|g'|_q=|g|_q$. Take a real number $p$ with $1/p+1/q=1$. Since $1<q\leq 2$, we have $2\leq p<\infty$. In particular, $q\leq p$. Thanks to H\"older's inequality,
\begin{align*}
|(f*g)(s)|&=\left|\sum_{\substack{|t|=k \\ |u|=\ell \\ s=tu}}f(t)g(u)\right| \\
&=\left|\sum_{|v|=j}f(t'v)g(v^{-1}u')\right| \\
&\leq\left(\sum_{|v|=j}|f(t'v)|^q\right)^{\frac{1}{q}}\left(\sum_{|v|=j}|g(v^{-1}u')|^p\right)^{\frac{1}{p}} \\
&\leq\left(\sum_{|v|=j}|f(t'v)|^q\right)^{\frac{1}{q}}\left(\sum_{|v|=j}|g(v^{-1}u')|^q\right)^{\frac{1}{q}} \\
&=f'(t')g'(u') \\
&=(f'*g')(s).
\end{align*}
Hence $|(f*g)\chi_m|\leq(f'*g')\chi_m$. Since $(k-j)+(\ell-j)=m$, it follows from the first part of the proof that
$$|(f*g)\chi_m|_q\leq|(f'*g')\chi_m|_q\leq|f'|_q|g'|_q=|f|_q|g|_q.$$
\epf

In the following lemma, the case where $p=q=2$ is given in the proof of {\cite[Theorem 1]{chh}}.

\blem\label{chh1}
Let $1\leq q\leq p\leq \infty$ with $1/p+1/q=1$. Let $\pi$ be a unitary representation of $\Gamma$ on a Hilbert space $\hil$ with a cyclic vector $\iota$ such that $\pi_{\xi,\xi}\in\ell_p(\Gamma)$. Then 
$$\|\pi(f)\|\leq\liminf_{n\to\infty} \left|(f^**f)^{(*2n)} \right|_q^{\frac{1}{4n}}$$
for $f\in c_c(\Gamma)$.
\elem

\bpf
Let $f\in c_c(\Gamma)$. In the proof of {\cite[Theorem 1]{chh}}, the following formula is obtained:
$$\|\pi(f)\|=\sup_{g\in c_c(\Gamma)}\lim_{n\to\infty}\left(\sum_{s\in\Gamma}(f^**f)^{(*2n)}(s)\langle\pi(s)\pi(g)\xi, \pi(g)\xi\rangle\right)^{\frac{1}{4n}}.$$
Fix  $g\in c_c(\Gamma)$ and we put $\varphi(s)=\langle\pi(s)\pi(g)\xi, \pi(g)\xi\rangle$. Note that 
$$\varphi(s)=\langle\pi(s)\pi(g)\xi, \pi(g)\xi\rangle=\sum_{t,u\in\free}\overline{g(u)}\pi_{\xi, \xi}(u^{-1}st)g(t)=(\overline{g}*\pi_{\xi, \xi}*g^\vee)(s),$$
where $g^\vee(s)=g(s^{-1})$. Consequently, $\pi_{\xi, \xi}\in\ell_p(\Gamma)$ implies $\varphi\in \ell_p(\Gamma)$. Then by H\"older's inequality, 
$$
\left|\sum_{s\in\Gamma}(f^**f)^{(2n)}(s)\varphi(s)\right|\leq\left|(f^**f)^{(2n)}\right|_q|\varphi|_p.
$$
Therefore it follows that
$$\|\pi(f)\|\leq\liminf_{n\to\infty}\left|(f^**f)^{(*2n)}\right|_q^{\frac{1}{4n}}.$$
\epf

By combining Lemma \ref{conv} and Lemma \ref{chh1}, we can prove the following.

\blem\label{rep}
Let $k$ be a non negative integer. Let $1\leq q\leq p\leq \infty$ with $1/p+1/q=1$. If a unitary representation $\pi$ of $\free$ on a Hilbert space $\hil$ has a cyclic vector $\xi$ such that $\pi_{\xi,\xi}\in \ell_p(\free)$, then 
$$\|\pi(f)\|\leq(k+1)|f|_q.$$
for $f\in c_c(\free)$ with $\mathrm{supp}(f)\subset W_k$.
\elem

\bpf
The case where $q=1$ and $p=\infty$ is trivial. So we may assume that $1<q\leq 2$ and $2\leq p<\infty$ with $1/p+1/q=1$. It is also shown by an argument similar as in {\cite[Lemma 1.4]{ha}}.

We consider the norm $\left|(f^**f)^{(*2n)}\right|_q$. Write $f_{2j-1}=f^*$ and $f_{2j}=f$ for $j=1, 2, \dots, 2n$. Then
$$(f^**f)^{(*2n)}=f_1*f_2*\cdots*f_{4n}.$$
We also denote $g=f_2*\cdots*f_{4n}$. So we have
$$(f^**f)^{(*2n)}=f_1*g.$$
Note that $\mathrm{supp}(f_j)\subset W_k$ for $j=1, 2, \dots, 4n$ and $g\in c_c(\free)$. Put $g_\ell=g\chi_\ell$. Then $\mathrm{supp}(g_\ell )\subset W_\ell$ and 
$$|g|_q^q=\sum_{\ell=0}^\infty|g_\ell|_q^q.$$
Here, remark that $|g_\ell|_q=0$ for all but finitely many $\ell$. Moreover set 
$$h=f_1*g=\sum_{\ell=0}^\infty f_1*g_\ell$$
and $h_m=h\chi_m$. Then $h\in c_c(\free)$ and
$$|h|_q^q=\sum_{m=0}^\infty|h_m|_q^q.$$
Here, notice that $|h_m|_q=0$ for all but finitely many $m$. By Lemma \ref{conv}, 
$$|(f_1*g_\ell)\chi_m|_q\leq|f_1|_q|g_\ell|_q$$
in the case where $|k-\ell|\leq m\leq k+\ell$ and $k+\ell-m$ is even. We also have 
$$|(f_1*g_\ell)\chi_m|_q=0$$
for any other values of $m$. Hence
\begin{align*}
|h_m|_q&=\left|\sum_{\ell=0}^\infty (f_1*g_\ell)\chi_m\right|_q \\
&\leq\sum_{\ell=0}^\infty \left|(f_1*g_\ell)\chi_m\right|_q \\
&=|f_1|_q\sum_{\substack{\ell=|m-k| \\ m+k-\ell \ \mathrm{even}}}^{m+k}|g_\ell|_q.
\end{align*}
By writing $\ell=m+k-2j$, 
\begin{align*}
|h_m|_q&\leq|f_1|_q\sum_{j=0}^{\min\{m, k\}}|g_{m+k-2j}|_q \\
&\leq|f_1|_q\left(\sum_{j=0}^{\min\{m, k\}}|g_{m+k-2j}|_q^q\right)^{\frac{1}{q}}\left(\sum_{j=0}^{\min\{m, k\}}1^p\right)^{\frac{1}{p}} \\
&\leq(k+1)^{\frac{1}{p}}|f_1|_q\left(\sum_{j=0}^{\min\{m, k\}}|g_{m+k-2j}|_q^q\right)^{\frac{1}{q}}.
\end{align*}
Then
\begin{align*}
|h|_q^q&=\sum_{m=0}^{\infty}|h_m|_q^q \\
&\leq(k+1)^{\frac{q}{p}}|f_1|_q^q\sum_{m=0}^\infty\sum_{j=0}^{\min\{m, k\}}|g_{m+k-2j}|_q^q \\
&=(k+1)^{\frac{q}{p}}|f_1|_q^q\sum_{j=0}^k\sum_{m=j}^{\infty}|g_{m+k-2j}|_q^q \\
&=(k+1)^{\frac{q}{p}}|f_1|_q^q\sum_{j=0}^k\sum_{\ell=k-j}^{\infty}|g_\ell|_q^q \\
&\leq(k+1)^{\frac{q}{p}}|f_1|_q^q\sum_{j=0}^k|g|_q^q \\
&=(k+1)^{\frac{q}{p}+1}|f_1|_q^q|g|_q^q.
\end{align*}
Hence $|f_1*g|_q\leq (k+1)|f_1|_q|g|_q$, i.e.,
$$|f_1*(f_2*\cdots*f_{4n})|_q\leq (k+1)|f_1|_q|f_2*\cdots*f_{4n}|_q.$$
Moreover, we inductively get
$$|(f^**f)^{(*2n)}|_q\leq(k+1)^{4n-1}|f|_q^{4n}.$$
Therefore it follows from Lemma \ref{chh1} that 
$$\|\pi(f)\|\leq\liminf_{n\to\infty}\left|(f^**f)^{(*2n)}\right|_q^{\frac{1}{4n}}\leq(k+1)|f|_q.$$
\epf

For a function $\varphi$ on $\Gamma$, we denote the corresponding linear functional on $c_c(\Gamma)$ by
$$\omega_\varphi(f)=\sum_{s\in\Gamma}f(s)\varphi(s)\ \ \ \mbox{for}\ f\in c_c(\Gamma).$$ 
Note that $\varphi$ is positive definite if and only if the functional $\omega_\varphi$ extends to a positive linear functional on $C^*(\Gamma)$, (see, e.g., \cite[Theorem 2.5.11]{boz}).

Now we can give a characterization of positive definite functions on $\free$, which can be extended to the positive linear functionals on $C_{\ell_p}^*(\free)$ for any $p\in[2,\infty)$. The case of $C_\lambda^*(\free)$ is given in {\cite[Theorem 3.1]{ha}}. We remind the reader that $C_\lambda^*(\free)= C_{\ell_2}^*(\free)$. Hence the following theorem is a generalization to the case of $C_{\ell_p}^*(\free)$ for any $p\in[2,\infty)$. 

For $0<\alpha<1$, we set $\varphi_\alpha(s)=\alpha^{|s|}$, and it is positive definite on $\free$ by \cite[Lemma 1.2]{ha}.

\bthm\label{posdef}
Let $2\leq p<\infty$. Let $\varphi$ be a positive definite function on $\free$. Then the following conditions are equivalent:
\benu
\item $\varphi$ can be extended to the positive linear functional on $C_{\ell_p}^*(\free)$.
\item $\sup_k|\varphi\chi_k|_p(k+1)^{-1}<\infty$.
\item The function $s\mapsto\varphi(s)(1+|s|)^{-1-\frac{2}{p}}$ belongs to $\ell_p(\free)$.
\item For any $\alpha\in(0,1)$, the function $s\mapsto\varphi(s)\alpha^{|s|}$ belongs to $\ell_p(\free)$.
\eenu
\ethm

\bpf
Without loss of generality, we may assume that $\varphi(e)=1$. The proof is based on the one in {\cite[Theorem 3.1]{ha}}.

(1)$\Rightarrow$(2): It follows from (1) that $\omega_\varphi$ extends to the state on $C_{\ell_p}^*(\free)$. Hence for $f\in c_c(\free)$, we have
$$\left|\omega_\varphi(f)\right|\leq\|f\|_{\ell_p}$$
Put 
$$f=|\varphi|^{p-2}\overline{\varphi}\chi_k.$$
Then 
$$\left|\omega_\varphi(f)\right|=|\varphi\chi_k|_p^p.$$
Notice that 
$$\|f\|_{\ell_p}=\sup\{\|\pi(f)\| : \pi\ \mbox{is an $\ell_p$-representation}\}.$$
Let $\pi$ be an $\ell_p$-representation of $\free$ on a Hilbert space $\hil$ with a dense subspace $\hil_0$. Then
$$\|\pi(f)\|^2=\sup_{\substack{\xi\in \hil_0 \\ \|\xi\|=1}}\langle\pi(f^**f)\xi, \xi\rangle_{\hil}.$$
Fix $\xi\in \hil_0$ with $\|\xi\|=1$. We denote by $\sigma$ the restriction of $\pi$ onto the subspace 
$$\hil_\sigma=\overline{\mathrm{span}}\{\pi(s)\xi : s\in\free\}\subset \hil.$$
Then
$$\langle\pi(f^**f)\xi, \xi\rangle_{\hil}=\langle\sigma(f^**f)\xi, \xi\rangle_{\hil_\sigma}.$$
Note that $\xi$ is cyclic for $\sigma$ such that $\sigma_{\xi, \xi}\in\ell_p(\free)$. Take a real number $q$ with $1/p+1/q=1$. Since $2\leq p<\infty$, we have $1<q\leq 2$. By Lemma \ref{rep}, 
$$\|\sigma(f)\|\leq(k+1)|f|_q.$$
Hence
$$\|\sigma(f^**f)\|=\|\sigma(f)\|^2\leq(k+1)^2|f|_q^2.$$
Therefore we obtain
$$\|f\|_{\ell_p}^2\leq(k+1)^2|f|_q^2=(k+1)^2|\varphi\chi_k|_p^{2(p-1)},$$
namely, 
$$\|f\|_{\ell_p}\leq(k+1)|\varphi\chi_k|_p^{p-1}.$$
Consequently,
$$|\varphi\chi_n|_p\leq k+1.$$

(2)$\Rightarrow$(3): 
\begin{align*}
\sum_{s\in\free}|\varphi(s)|^p(1+|s|)^{-p-2}&=\sum_{k=0}^\infty\sum_{|s|=k}|\varphi(s)|^p(1+k)^{-p-2} \\
&=\sum_{k=0}^\infty|\varphi\chi_k|_p^p(1+k)^{-p}(1+k)^{-2} \\
&\leq\left\{\sup_k|\varphi\chi_k|_p(k+1)^{-1}\right\}^p\sum_{k=0}^\infty\frac{1}{(k+1)^2}<\infty.
\end{align*}

(3)$\Rightarrow$(4): Easy.

(4)$\Rightarrow$(1): 
Note that $\psi_\alpha(s)=\varphi(s)\alpha^{|s|}$ is also positive definite. By the GNS construction, we obtain the unitary representation $\sigma_\alpha$ of $\free$ with the cyclic vector $\xi_\alpha$ such that
$$\psi_\alpha(s)=\langle\sigma_\alpha(s)\xi_\alpha, \xi_\alpha\rangle.$$
Since $\sigma_\alpha$ is an $\ell_p$-representation, $\psi_\alpha$ can be seen as a state on $C_{\ell_p}^*(\free)$. By taking the weak-$*$ limit of $\psi_\alpha$ as $\alpha\nearrow 1$, we conclude that $\varphi$ can be extended to the state on $C_{\ell_p}^*(\free)$.
\epf

\bcor\label{exp}
Let $p\in[2, \infty)$ and $\alpha\in(0, 1)$. The positive definite function $\varphi_\alpha$ can be extended to the state on $C_{\ell_p}^*(\free)$ if and only if 
$$\alpha\leq(2d-1)^{-\frac{1}{p}}.$$
\ecor

\bpf
Note that 
\begin{eqnarray*}
\varphi_\alpha\in\ell_p(\free)&\Longleftrightarrow& \sum_{k=1}^\infty(2d-1)^{k-1}\alpha^{pk}<\infty \\
&\Longleftrightarrow& (2d-1)\alpha^p<1 \\
&\Longleftrightarrow& \alpha<(2d-1)^{-\frac{1}{p}}.
\end{eqnarray*}
Hence the corollary follows from Theorem \ref{posdef}.
\epf

\brem
Let $\pi_\alpha$ be the GNS representation of $\varphi_\alpha$. Then $\pi_\alpha$ is weakly contained in $\lambda$ if and only if $\alpha\leq(2d-1)^{-\frac{1}{2}}$. For $\alpha>(2d-1)^{-\frac{1}{2}}$, we refer the reader to \cite{ps} and \cite{sz}.
\erem

As a consequence, we can obtain the following result. See also \cite[Proposition 4.4]{bg}.

\bcor\label{non-iso}
For $2\leq q<p\leq\infty$, the canonical quotient map from $C_{\ell_p}^*(\free)$ onto $C_{\ell_q}^*(\free)$ is not injective.
\ecor

\bpf
It suffices to consider the case where $p\ne\infty$, because $\free$ is not amenable. 

Suppose that the canonical quotient map from $C_{\ell_p}^*(\free)$ onto $C_{\ell_q}^*(\free)$ is injective for some $q<p$. Take a real number $\alpha$ with 
$$(2d-1)^{-\frac{1}{q}}<\alpha\leq(2d-1)^{-\frac{1}{p}}.$$ 
By using Corollary \ref{exp}, 
$$\left|\omega_{\varphi_\alpha}(f)\right|\leq\|f\|_{\ell_p}=\|f\|_{\ell_q}\ \ \ \mbox{for}\ f\in c_c(\free).$$
Therefore it follows that $\varphi_\alpha$ can be also extended to the state on $C_{\ell_q}^*(\free)$, but it contradicts to the choice of $\alpha$. 
\epf

\brem
The previous result has also shown by N. Higson and N. Ozawa, independently. See also \cite[Remark 4.5]{bg}.
\erem

In \cite{p}, Powers proves that $C_\lambda^*(\mathbb{F}_2)$ has a unique tracial state. In \cite{ha}, Haagerup gives another proof of the uniqueness. Thanks to Theorem \ref{posdef}, Haagerup's argument also works for the case of $C_{\ell_p}^*(\free)$.

\bcor\label{trace}
For each $p\in[2,\infty)$, the $C^*$-algebra $C_{\ell_p}^*(\free)$ has a unique tracial state $\tau\circ\lambda_p$, where $\lambda_p$ is the canonical quotient map $C^*_\lp(\free)$ onto $C^*_\lambda(\free)$ and $\tau$ is the unique tracial state on $C^*_\lambda(\free)$.
\ecor

\bpf
Any tracial state on $C_{\ell_p}^*(\free)$ corresponds a positive definite function on $\free$. Take such a positive definite function $\varphi$. Then $\varphi|_K$ is constant for any conjugacy class $K$ in $\free$. Take a conjugacy class $K$ in $\free$ such that $K\ne\{e\}$. Put $k=\min\{|s| : s\in K\}$. Then $|W_{k+2n}\cap K|\geq (2d-1)^{n-1}$ for $n\geq 1$. Hence if $\varphi|_K=C$ for some non-zero constant $C$, then
$$
\sum_{s\in C}|\varphi(s)|^p\alpha^{p|s|}=\sum_{n=0}^\infty\sum_{s\in W_{k+2n}\cap K}C^p\alpha^{pn}\geq C^p\sum_{n=0}^\infty(2d-1)^{n-1}\alpha^{pn}=\infty
$$
for $\alpha\geq(2d-1)^{-\frac{1}{p}}$. This contradicts (4) in Theorem \ref{posdef}.
\epf

We define two algebraic ideals by
$$D_p^-(\Gamma)=\bigcup_{\varepsilon>0}\ell_{p-\varepsilon}(\Gamma)$$
for $1<p\leq\infty$, and
$$D_p^+(\Gamma)=\bigcap_{\varepsilon>0}\ell_{p+\varepsilon}(\Gamma)$$
for $1\leq p<\infty$. Note that $\ell_q(\Gamma)\subsetneqq D_p^-(\Gamma)\subsetneqq\ell_p(\Gamma)$ for $1\leq q<p\leq\infty$, and $\ell_p(\Gamma)\subsetneqq D_p^+(\Gamma)\subsetneqq\ell_q(\Gamma)$ for $1\leq p<q\leq\infty$. Then we also obtain the following.

\bcor
\benu
\item For $2\leq p<\infty$, the $C^*$-algebra $C_\lp^*(\free)$ is canonically isomorphic to $C_{D_p^+}^*(\free)$. In particular, $C_\lambda^*(\free)=C_{D_2^+}^*(\free)$.
\item For $2< p\leq\infty$, the $C^*$-algebra $C_\lp^*(\free)$ is canonically isomorphic to $C_{D_p^-}^*(\free)$. In particular, $C^*(\free)=C_{D_\infty^-}^*(\free)$.
\eenu
\ecor

\bpf
(1) It suffices to show that if $\varphi$ is a positive definite function on $\free$ which can be extended the positive linear functional  on $C_{D_p^+}^*(\free)$, then $\varphi$ can be also extended to the one on $C_\lp^*(\free)$. 

Now assume that $\varphi$ is a positive definite function on $\free$, which can be extend to the positive linear functional on $C_{D_p^+}(\free)$. Then for any $q\in(p, \infty)$, $\varphi$ can be also extended to the positive linear functional on  $C_{\ell_q}^*(\free)$. By Theorem \ref{posdef}, $\varphi\varphi_\alpha\in\ell_q(\free)$ for any $\alpha\in(0, 1)$. We set 
$$r=\frac{pq}{q-p}.$$
Then $1/p=1/q+1/r$. If we take $(2d-1)^{-2/r}<\beta<(2d-1)^{-1/r}$, then $\varphi_\beta\in\ell_r(\free)$. Since
$$|\varphi\varphi_\alpha\varphi_\beta|_p\leq|\varphi\varphi_\alpha|_q|\varphi_\beta|_r,$$
we have $\varphi\varphi_\alpha\varphi_\beta\in\ell_p(\free)$. Namely $(\varphi\varphi_\beta)\varphi_\alpha\in\ell_p(\free)$ for any $\alpha\in(0,1)$. Thus $\varphi\varphi_\beta$ can be extended to the positive linear functional on $C_\lp^*(\free)$. If $q\searrow p$, then $r\nearrow \infty$ and $\beta\nearrow 1$. Hence $\varphi\varphi_\beta\to\varphi$ in the weak-$*$ topology. Therefore $\varphi$ can be also extended to the positive linear functional on $C_\lp^*(\free)$.

(2) The proof is quite similar as in (1). Assume that $\varphi$ is a positive definite function on $\free$ which can be extended to the positive linear functional on $C_\lp^*(\free)$. By Theorem \ref{posdef}, we have $\varphi\varphi_\alpha\in\ell_p(\free)$ for any $\alpha\in(0, 1)$. For any $q\in[2, p)$, we set 
$$r=\frac{pq}{p-q}.$$
Then $1/q=1/p+1/r$. If we take $(2d-1)^{-2/r}<\beta<(2d-1)^{-1/r}$, then $\varphi_\beta\in\ell_r(\free)$. Since 
$$|\varphi\varphi_\alpha\varphi_\beta|_q\leq|\varphi\varphi_\alpha|_p|\varphi_\beta|_r,$$
we have $\varphi\varphi_\alpha\varphi_\beta\in\ell_q(\free)$. Namely, $(\varphi\varphi_\beta)\varphi_\alpha\in\ell_q(\free)$ for any $\alpha\in(0,1)$. Thus $\varphi\varphi_\beta$ can be extended to the positive linear functional on $C_{\ell_q}^*(\free)$, and so $\varphi\varphi_\beta$ can be extended to the one on $C_{D_p^-}^*(\free)$. If $q\nearrow p$, then $r\nearrow \infty$ and $\beta\nearrow 1$. Hence $\varphi\varphi_\beta\to\varphi$ in the weak-$*$ topology. Therefore $\varphi$ can be also extended to the positive linear functional on $C_{D_p^-}^*(\free)$.
\epf

\brem
The isomorphism $C_\lambda^*(\Gamma)=C_{D_2^+}(\Gamma)$ has been already obtained in \cite{bg} for any $\Gamma$. 
\erem

\bibliographystyle{amsplain}

\end{document}